\documentclass[11pt,a4paper]{amsart}
\usepackage[utf8]{inputenc}
\usepackage[english]{babel}
\usepackage{amsthm}
\usepackage{amssymb}
\usepackage{enumerate}

\newcommand{\N}{\mathbb{N}}
\newcommand{\Z}{\mathbb{Z}}

\newcommand{\End}{\text{End}}
\newcommand{\Rat}{\text{Rat}}
\newcommand{\Orb}{\text{Orb}}
\newcommand{\Rec}{\text{Rec}}
\newcommand{\Alg}{\text{Alg}}
\newcommand{\CF}{\text{CF}}
\newcommand{\Aut}{\text{Aut}}

\newcommand{\mc}{\mathcal}

\theoremstyle{definition}
\newtheorem{theorem}{Theorem}[section]

\newtheorem{proposition}[theorem]{Proposition}

\newtheorem{problem}[theorem]{Problem}
\newtheorem{lemma}[theorem]{Lemma}

\newtheorem{remark}[theorem]{Remark}

\begin{document}
 
 
\title{Subsets of groups in public-key cryptography}
\author{André Carvalho}
\address[A. Carvalho]{%
Center for Mathematics and Applications (NOVA Math)\\
NOVA School of Science and Technology\\
NOVA University of Lisbon\\
2829--516 Caparica\\
Portugal
}
\email{%
andrecruzcarvalho@gmail.com
}
\author{António Malheiro}
\address[A. Malheiro]{%
Center for Mathematics and Applications (NOVA Math) \& Department of Mathematics\\
NOVA School of Science and Technology\\
NOVA University of Lisbon\\
2829--516 Caparica\\
Portugal
}
\email{%
ajm@fct.unl.pt
}
\thanks{This work is funded by national funds through the FCT - Fundação para a Ciência e a Tecnologia, I.P., under the
  scope of the projects UIDB/00297/2020 and UIDP/00297/2020 (Center for Mathematics and Applications).}

\begin{abstract}
We suggest the usage of algebraic subsets instead of subgroups in public-key cryptography. In particular, we present the subset version of two protocols introduced by Shpilrain and Ushakov with some examples in ascending HNN-extensions of free-abelian groups and discuss their resistance to length and distance based attacks. We  also introduce several new group theoretic problems arising from this work.
\end{abstract}

\maketitle

\section{Introduction}
In recent years, there has been an increasing interest in the development of group-based protocols for public-key cryptography with the appearance of  different types of protocols, such as key exchange and authentication protocols and zero-knowledge proofs (see, for example, \cite{[FKN23]} for a survey). We will focus on key exchange protocols, but, in theory, our approach serves other cryptosystems in the same way. Typically, when defining a key exchange protocol, there is a (difficult) algorithmic problem about the group involved in the construction of the key: in \cite{[AAG99]}, Anshel, Anshel and Goldfeld proposed a key exchange protocol on braid groups, involving a problem resembling the subgroup-restricted simultaneous conjugacy problem; in  \cite{[KLCHKP00]}, Ko, Lee et al. suggest a protocol similar to Diffie--Hellman's  based on the conjugacy problem in braid groups; Shpilrain and Ushakov proposed in \cite{[SU05],[SU06]} two protocols based on the decomposition problem for the Thompson's group $F$ and  braid groups; in \cite{[SZ06]}, Shpilrain and Zapata propose a cryptosystem based on the membership problem in free metabelian groups; and there are many others  (for a general reference we refer the reader to \cite{[MSU11],[FKN23]}). In most cases, the decision problem is not where the security of the protocol relies, but on the \emph{search} version of the problem, where we look for a witness knowing that a given property holds.
For the purpose of this paper, the most relevant protocols will be the ones  in \cite{[SU05],[SU06]} based on the decomposition search problem. 

The study of subgroups is as old as group theory. The study of subsets of groups, however, is much more recent. Naturally, when considering arbitrary subsets of a group not much can be said, as there is no structure. However, when the subset is defined via some language theoretic condition, interesting results can be obtained (we refer the reader to \cite{[Ber79],[BS21]}). Variations of decision problems involving rational and algebraic subsets are in general harder (at least not easier) than the ones concerning finitely generated subgroups, so it is natural to consider the search versions of these problems as well, and propose them as possible security assumptions of certain cryptosystems. More details are presented in Section \ref{sec: preliminaries}. 

In this paper, we consider variants of the protocols   in \cite{[SU05],[SU06]} replacing finitely generated subgroups by subsets of the group. We show that using rational subsets does not make the protocols more secure, but algebraic subsets (the context-free counterpart of rational subsets) appear to do so. Implementation can be easily done via pushdown automata or context-free grammars.
 We present the subsets variants of both protocols, provide some examples and an heuristic argument as to why these variations should be more resistant to traditional length and distance based attacks. 

In Section \ref{sec: preliminaries}, we introduce some definitions and known results on subsets of groups and on algorithmic problems (decision and search versions). In Section \ref{protocols}, we explicitly describe the protocols obtained by modifying Shpilrain and Ushakov's protocols from  \cite{[SU05],[SU06]} by replacing finitely generated subgroups with algebraic subsets. In Section \ref{sec: other protocols} we discuss the application of this method in other protocols showing other examples  for which it works and demonstrating an obstruction to the usage of subsets in Anshel--Anshel--Goldfeld's protocol. In Section \ref{sec: choosing}, we give some key aspects to consider when choosing the subsets and in Section \ref{proposal}, we present a protocol considering algebraic subsets of ascending HNN-extensions of free-abelian groups. In Section \ref{sec: attacks} we present an heuristic argument to emphasize the resistance of this method against length and distance based attacks. In Section \ref{sec: other problems}, we present some interesting orbit problems arising from this work and we end with some final remarks in Section \ref{sec: conclusion}.

\section{Preliminaries}\label{sec: preliminaries}
In this section we review some concepts on subsets of groups and on algorithmic problems in group theory. 
\subsection{Subsets of groups}

Let $G=\langle A\rangle$ be a finitely generated group, $A$ be a finite generating set, $\tilde A=A\cup A^{-1}$ and $\pi:\tilde A^*\to G$ be the canonical (surjective) homomorphism. This notation will be kept throughout the paper.

A subset $K\subseteq G$ is said to be \emph{rational} if there is some rational language $L\subseteq \tilde A^*$ such that $L\pi=K$ and \emph{recognizable} if $K\pi^{-1}$ is rational.

We will denote by $Rat(G)$ and $Rec(G)$ the class of rational and recognizable subsets of $G$, respectively. 
Rational subsets generalize the notion of finitely generated subgroups.

\begin{theorem}[\cite{[Ber79]}, Theorem III.2.7]
\label{AnisimovSeifert}
Let $H$ be a subgroup of a group $G$. Then $H\in \text{Rat } G$ if and only if $H$ is finitely generated.
\end{theorem}

Similarly, recognizable subsets generalize the notion of finite index subgroups.

\begin{proposition}[\cite{[Ber79]}, Exercise III.1.3]
\label{rec fi}
Let $H$ be a subgroup of a group $G$. Then $H\in \text{Rec } G$ if and only if $H$ has finite index in $G$.
\end{proposition}

In fact, it is easy to see that if $G$ is a group and $K$ is a subset of $G$ then $K$ is recognizable if and only if $K$ is a (finite) union of cosets of a subgroup of finite index.

A natural generalization concerns the class of context-free languages. 
A subset $K\subseteq G$ is said to be \emph{algebraic} if there is some context-free language $L\subseteq \tilde A^*$ such that $L\pi=K$ and \emph{context-free} if $K\pi^{-1}$ is context-free. 
We will denote by $\Alg(G)$ and $\CF(G)$ the class of algebraic and context-free subsets of $G$, respectively. 
 It follows from \cite[Lemma 2.1]{[Her91]} that $\CF(G)$ and $\Alg(G)$ do not depend on the alphabet $A$ or the surjective homomorphism $\pi$.

For a finitely generated group $G$, it is immediate from the definitions that $\Rec(G)\subseteq \CF(G) \subseteq \Alg(G)$ and that $\Rec(G)\subseteq \Rat(G) \subseteq \Alg(G).$ It is proved in \cite{[Her91]} that 
$$\CF(G)=\Alg(G) \Leftrightarrow \CF(G)=\Rat(G) \Leftrightarrow \text{ G is virtually cyclic.}$$

However, there is no general inclusion between $\Rat(G)$ and $\CF(G)$. For example, if $G$ is virtually abelian, then $\CF(G)\subseteq \Alg(G)= \Rat(G)$ (and the inclusion is strict if the group is not virtually cyclic) and if the group is virtually free, then $\Rat(G)\subseteq \CF(G)$ (see \cite{[Her91]}).

We end the subsection with a useful lemma.
\begin{lemma}\label{lem: generated algebraic}
Let $G$ be a group and $S\in \Alg(G)$. Then $\langle S\rangle\in \Alg(G).$
\end{lemma}
\noindent\textit{Proof.} Let $A$ be a generating set for $G$ and $L\subseteq \tilde A^*$ be a context-free language such that $L\pi=S$. Then, $(L\cup L^{-1})^*\subseteq \tilde A^*$ is a context-free language and  $(L\cup L^{-1})^*\pi=\langle S\rangle$, so $\langle S\rangle\in \Alg(G).$
\qed

\subsection{Decision vs search problems}
Decision problems play a big role in group theory. These consist on deciding whether a given property holds for a certain input. For example, we may want to decide whether two elements are conjugate or two subgroups isomorphic. More recently, variants of many standard decision problems involving subsets of groups have been studied. For example, instead of the standard membership problem (involving subgroups), we may study the rational subset membership problem (see \cite{[Loh15]} for a survey). This modification is relevant: for instance, the rational membership problem  is undecidable for the free metabelian group with two generators and for the right-angled Artin group $\mathbb G(P_4)$, while the subgroup membership is decidable for these groups \cite{[LS08]}. Variations of the conjugacy problem have also been considered in \cite{[LS11],[Car23b]} in two different ways: we may want to decide whether an element has a conjugate in a given subset or whether two elements are conjugate with a conjugator belonging to a certain subset. Similar variants have been considered for the intersection, twisted conjugacy or equality problems, both for groups and semigroups \cite{[AH89],[Sil23],[Car23b]}. Since rational subsets generalize the notion of finitely generated subgroups, decision problems in rational or algebraic subsets are in general harder than their subgroups counterpart.

In group-based cryptography, it is common that the security of a protocol relies on a \emph{search} version of a decision problem. This consists on, knowing that a certain property holds for our  input, searching for a witness. In most cases, these problems are decidable, what we care about is the complexity. Search problems have not been studied for long by group theorists and the interest in them grew because of cryptography. Therefore, it seems to make sense to consider the study of the subsets-variation of search problems  involving subgroups.
\begin{problem}
Study the complexity of the (rational, algebraic) subset restricted conjugacy search problem, factorization problem and decomposition problem for several classes of groups.
\end{problem}

The rational subsets version of these search problems is not easier than the finitely generated subgroup version, as it is not hard to construct an automaton that allows us to view a finitely generated subgroup as a rational subset. While there is not much work done on the factorization and decomposition problems, the subset-restricted conjugacy decision problem is highly nontrivial: for finitely generated virtually free groups it is known to be decidable (see \cite{[LS11]}), but it is not known to be decidable in free-by-cyclic groups.

\section{Shpirlain--Ushakov's Protocols} \label{protocols}

In \cite{[SU05]}, Shpilrain and Ushakov propose the following key exchange system and they suggest the Thompson's group $F$ as the platform group to be used. Matucci proved in \cite{[Mat08]} that this setting could be broken deterministically and Ruinskiy--Shamir--Tsaban use distance functions to break the protocol (non-deterministically) in \cite{[RST07]}. The protocol is as follows:

\begin{center}
\fbox{\begin{minipage}{\textwidth}
\noindent\textbf{Public information:} group $G$, $w\in G$, $A,B\leq_{f.g.} G$ such that $ab=ba$ for all $(a,b)\in A\times B$.\\

\noindent\textbf{Key exchange:}
\begin{enumerate}
\item Alice picks $a_1\in A$, $b_1\in B$ and sends $a_1wb_1$ to Bob.
\item Bob picks  $a_2\in A$, $b_2\in B$ and sends $b_2wa_2$ to Alice.
\item Alice computes $K_A=a_1(b_2wa_2)b_1$ and Bob computes $K_B=b_2(a_1wb_1)a_2$.\\
\end{enumerate}

\noindent\textbf{Shared secret key:}  $K=K_A=K_B=a_1b_2wa_2b_1$\\
\end{minipage}}
\end{center}

The security of this protocol seems to be based on the \emph{subgroup-restricted decomposition search problem}:\\

Given two elements $w,w'\in G$ and two subgroups $A,B\leq G$, find two elements $x\in A$ and $y\in B$ such that $xwy=w'$, provided at least one such pair of elements exists.\\

As noted in \cite{[MSU11]}, this problem is equivalent to the  \emph{subgroup-restricted factorization search problem}:\\

Given an element $w\in G$ and two subgroups $A,B\leq G$, find two elements $x\in A$ and $y\in B$ such that $xy=w$, provided at least one such pair of elements exists.\\

That is the case, since there are elements $x\in A$ and $y\in B$ such that $xwy=w'$ if and only if there are elements $x\in A^w$ and $y\in B$ such that $xy=w^{-1}w'$.

\begin{remark}\label{rmk: centralizer}
In order to break the protocol, an adversary Eve does not have to necessarily obtain $a_1',a_2'\in A$ and $b_1',b_2'\in B$ such that $a_1'wb_1'=a_1wb_1$ and $b_2'wa_2'=b_2wa_2$, but \emph{only} $a\in C(B)$, $b\in C(A)$ and $c,d\in G$ such that $awb=a_1wb_1$ and $cwd=b_2wa_2$. Indeed, if that is the case, then 
$$acwdb=ab_2wa_2b=b_2awba_2=b_2a_1wb_1a_2.$$
Of course, knowing there are solutions in the subgroups $A$ and $B$ usually allows us to make the search more effective.
\end{remark}

We propose a variation on this scheme that can theoretically be applied in many other schemes involving subgroups that consists on replacing the subgroups $A$ and $B$ by subsets. Naturally, we want these subsets to be infinite but definable by a finite amount of information. 
Rational subsets are  the most standard and the most well-studied in the literature, so a natural candidate would be replacing finitely generated subgroups by rational subsets in the above protocol.

 However, it follows from Remark \ref{rmk: centralizer}, that this is not  harder to break than the previous version of the protocol.
 Indeed, suppose that we have rational subsets $A,B\in \Rat(G)$ such that $[A,B]=1$. Then, by \cite[Theorem 4.2]{[Gil95]}, $\langle A \rangle$ and  $\langle B \rangle$ are finitely generated and  finite sets of generators for these subgroups are computable. Also, we have that $C(A)=C(\langle A\rangle)$ and $C(B )=C(\langle B\rangle)$, and so, if Eve finds   $a_1',a_2'\in \langle A\rangle$ and $b_1',b_2'\in \langle B \rangle$ such that $a_1'wb_1'=a_1wb_1$ and $b_2'wa_2'=b_2wa_2$, then 
 \begin{equation*}
 a_1'b_2'wa_2'b_1'=a_1'b_2wa_2b_1'=b_2a_1'wb_1'a_2=b_2a_1wb_1a_2.
 \end{equation*}
 Notice that the second equality follows from the fact that $b_2\in C(A)$, and so $b_2\in C(\langle A \rangle)$.

 Algebraic subsets can have a much more complicated structure and are also not hard to define by a finite structure such as  a context-free grammar or a pushdown automaton and are, in general, significantly more complex than rational subsets. 

It is proved in \cite{[Car23c]} that although it is true that, for $H\leq G$, $$\Rat(H)=\{K\subseteq H \mid K\in \Rat(G)\},$$ the same property does not hold for algebraic subsets. This means that algebraic subsets of a supergroup can be more complicated than algebraic subsets of a subgroup, while that does not happen in the rational case. Also, it is observed in \cite{[Car23c]} that, given any group $G$, any element $g\in G$ and any automorphism $\phi\in \Aut(G)$,  the orbit of $g$ under the iteration of $\phi$ is an algebraic subset of $G\rtimes_\phi \Z$, via the language $\{t^{-n}wt^n\mid n\in \Z\}$, where $w$ is a word representing $g$. So, being able to decide membership in algebraic subsets of a $G$-by-$\Z$ group implies decidability of Brinkmann's equality problem on $G$, which is a difficult problem even in the free group \cite{[Bri10]}, and so the factorization decision problem relative to algebraic subsets of $G$-by-$\Z$ groups should also be difficult (in the case where $B$ is trivial, this is simply deciding membership in $A$).

So, the protocol would now be the following:
\begin{center}
\fbox{\begin{minipage}{\textwidth}
\noindent\textbf{Public information:} group $G$, $w\in G$, context free grammars $\mc G_1$ and $\mc G_2$ such that $ab=ba$ for all $(a,b)\in L(\mc G_1)\pi\times L(\mc G_2)\pi$. \\

\noindent\textbf{Key exchange:}
\begin{enumerate}
\item Alice picks $a_1\in L(\mc G_1)$, $b_1\in L(\mc G_2)$ and sends $a_1wb_1$ to Bob.
\item Bob picks  $a_2\in L(\mc G_1)$, $b_2\in L(\mc G_2)$ and sends $b_2wa_2$ to Alice.
\item Alice computes $K_A=a_1(b_2wa_2)b_1$ and Bob computes $K_B=b_2(a_1wb_1)a_2$. \\
\end{enumerate}

\noindent\textbf{Shared secret key:}  $K=K_A=K_B=a_1b_2wa_2b_1$\\
\end{minipage}}
\end{center}

In \cite{[SU06]}, Shpilrain and Ushakov propose a variation of the above protocol, which can be seen as a version where one of the subgroups is hidden:

\begin{center}
\fbox{\begin{minipage}{\textwidth}
\noindent\textbf{Public information:}  group $G$, $w\in G$.\\

\noindent\textbf{Key exchange:}
\begin{enumerate}
\item Alice picks $a_1\in G$, chooses a subgroup of $C_G(a_1)$ and publishes its generators $A=\{\alpha_1,\ldots,\alpha_k\}$
\item Bob picks  $b_2\in G$,  chooses a subgroup of $C_G(b_1)$ and publishes its generators $B=\{\beta_1,\ldots,\beta_k\}$
\item Alice chooses $a_2\in \langle\beta_1,\ldots,\beta_k \rangle$ and sends $a_1wa_2$ to Bob
\item Bob chooses $b_1\in \langle\alpha_1,\ldots,\alpha_k \rangle$ and sends $b_1wb_2$ to Alice
\item Alice computes $K_A=a_1(b_1wb_2)a_2$ and Bob computes $K_B=b_1(a_1wa_2)b_2$. 
\end{enumerate}

\noindent\textbf{Shared secret key:}   $K=K_A=K_B=a_1b_1wa_2b_2$\\
\end{minipage}}
\end{center}

The key points is that $a_1$ and $a_2$ are chosen arbitrarily, so if an adversary wants to restrict the search to a given subgroup, she would have to compute $C_G(B)$ and $C_G(A)$, since $a_1\in C_G(B)$ and $b_1\in C_G(A)$. This is in general hard to compute since these are the intersection of centralizers of elements $b_i$ and $a_i$. Also, such subgroups might not be finitely generated: for example, among all solvable groups, polycyclic groups are exactly those for which the centralizer of any finitely generated subgroup is finitely generated \cite{[Len96]}. This means that, for instance in the Baumslag-Solitar group $BS(1,2)$, there are finitely generated subgroups whose centralizer is not finitely generated. Moreover, $C_G(a_1)$ might not be finitely generated itself: for example  the Grigorchuk group has generators with a non-finitely generated centralizer \cite{[Roz94]}. The fact that $C_G(a_1)$ might not be finitely generated itself is not an obstruction to this method, as Alice only needs to compute a finite set of elements in the centralizer (not necessarily a generating set).

The subset version would be the following:

\begin{center}
\fbox{\begin{minipage}{\textwidth}
\noindent\textbf{Public information:}  group $G$, $w\in G$.\\

\noindent\textbf{Key exchange:}
\begin{enumerate}
\item Alice picks $a_1\in G$, chooses an algebraic subset of $C_G(a_1)$ and publishes a CFG $\mc G$ recognizing it
\item Bob picks  $b_2\in G$,  chooses an algebraic subset of $C_G(b_2)$ and publishes a CFG $\mc H$ recognizing it
\item Alice chooses $a_2\in L(\mc H)\pi$ and sends $a_1wa_2$ to Bob
\item Bob chooses $b_1\in L(\mc G)\pi$ and sends $b_1wb_2$ to Alice
\item Alice computes $K_A=a_1(b_1wb_2)a_2$ and Bob computes $K_B=b_1(a_1wa_2)b_2$.
\end{enumerate}

\noindent\textbf{Shared secret key:}   $K=K_A=K_B=a_1b_1wa_2b_2$\\
\end{minipage}}
\end{center}

In order to generate the CFG, Alice and Bob could proceed in the same way as done in the subgroup case. More details on the generation of a good algebraic subset will be given in Section \ref{sec: choosing}. However, computing the centralizer of an algebraic subset should be a very difficult task, much more complicated than the finitely generated subgroup case. The idea is that, for finitely generated subgroups, this is the intersection of centralizers of elements, which, despite being difficult, can be done for some classes of groups. However, the authors are not aware of nontrivial results in this direction for any class of groups regarding the computation of centralizers of algebraic subsets. 
\begin{problem}
Find examples of non-virtually abelian groups for which centralizers of algebraic subsets are computable.
\end{problem}
If the subsets are well chosen, this should represent a significant obstacle to an adversary and an increment in security when compared to the subgroup version.

\section{Application to other protocols}\label{sec: other protocols}
This idea could in theory be applied to any protocol involving finitely generated subgroups. For example,  the protocols in \cite[Section 4.3]{[MSU11]} and \cite{[KLCHKP00],[Kur06]} and the authentication protocol in \cite{[LC05]} should be generalizable in the same way. 

However, one must be cautious when doing so. We now give examples of a protocol for which this approach  cannot be easily adapted.

Anshel--Anshel--Goldfeld's (AAG) protocol was introduced in \cite{[AAG99]} and can be described as follows:

\begin{center}
\fbox{\begin{minipage}{\textwidth}
\noindent\textbf{Public information:}  group $G$, elements $a_1,\ldots, a_k, b_1,\ldots, b_m\in G$.\\

\noindent\textbf{Key exchange:}
\begin{enumerate}
\item Alice picks $x\in G$ as a word in $a_1,\ldots, a_k$ and  sends $b_1^x,\ldots, b_m^x$ to Bob.
\item Bob picks $y\in G$ as a word in $b_1,\ldots, b_m$ and  sends $a_1^y,\ldots, a_k^y$ to Alice.
\item Alice computes $x(a_1^y,\ldots, a_k^y)=x^y$ and Bob computes $y(b_1^x,\ldots, b_m^x)=y^x$. 
\end{enumerate}

\noindent\textbf{Shared secret key:}   $K=x^{-1}y^{-1}xy$ (this can be computed by Bob and Alice).\\
\end{minipage}}
\end{center}

Taking subsets instead of the generators $a_1,\ldots, a_k$ and $b_1,\ldots, b_m$ does not help. If Alice picks $x$ as a word in a context-free language, it is not clear what she should send Bob. If we fix subsets $A,B\subseteq G$ and we make Alice pick a word in a context-free language $L\subseteq \tilde A^*$  and Bob picks $y$ as a word in a context-free language $K\subseteq \tilde B^*$, if the adversary Eve manages to obtain $x'=_Gx$ and $y'=_Gy$ as words in $\tilde A$ and $\tilde B$, respectively, then she can compute $x'(a_1^{y'},\ldots, a_k^{y'})={x'}^{y'}=x^y$ and obtain the shared secret key, reducing the problem to the one in $\langle A\rangle$ and $\langle B\rangle$.

\section{Choosing the subsets}
\label{sec: choosing}

\subsection{General assumptions}
We remark that the obstruction to the usage of rational subsets does not apply to the algebraic setting. Indeed, letting $F_2$ denote the free group of rank 2, we have that $S=\{a^nba^{-n}\mid n\in \Z\}\in \Alg(F_2)$ and $\langle S\rangle$ is not finitely generated. 

\subsubsection{First protocol}
Let $S,T$ be our public algebraic subsets. We have that $S\subseteq C(T)$. In order to make this protocol stronger than the version concerning finitely generated subgroups, it is important to make it hard for an adversary to find  finitely generated subgroups $H_1$  and $H_2$ such that $S\subseteq H_1\subseteq C(T)$ and $T\subseteq H_2\subseteq C(S)$. In particular $\langle S \rangle$ and $\langle T \rangle$  should not be finitely generated. 
Indeed, if such  subgroups $H_1$ and $H_2$ can be computed, then if the adversary Eve can break the subgroup version of the protocol, by finding $a_1', a_2'\in H_1$ and $b_1',b_2'\in H_2$ such that 
\begin{align*}
a_1wb_1=a_1'wb_1' \quad\text{ and  }\quad b_2wa_2=b_2'wa_2',
\end{align*}
 which necessarily exist since $S\subseteq H_1$ and $T\subseteq H_2$, then $$a_1'b_2'wa_2'b_1'=a_1'b_2wa_2b_1'=b_2a_1'wb_1'a_2=b_2a_1wb_1a_2$$
 and the shared private key could be computed. Notice that $a_1'b_2=b_2a_1'$ and $a_2b_1'=b_1'a_2$ since $H_1\subseteq C(T)$ and $H_2\subseteq C(S)$.

We call this the \emph{Middle Subgroup Problem (MSP)}: \\

Given algebraic subsets $S,T\in \Alg(G)$, can we decide if there is some finitely generated subgroup $H\leq G$ such that $S\subseteq H\subseteq C(T)$?\\

Naturally, we can define the  \emph{Middle Subgroup Search Problem (MSSP)}:\\

Given algebraic subsets $S,T\in \Alg(G)$, can we compute words $h_1,\ldots, h_k$ such that  $S\subseteq \langle h_1,\ldots, h_k\rangle\subseteq C(T)$, provided that such words exist?\\
 
 This condition is important to ensure that, even in the case where the subgroup version of the protocol can be broken for \emph{every} choice of subgroups, the algebraic subsets version might not  be. The ideal setting is then the case where the MSP outputs \texttt{NO}, but, of course it might be the case where  subgroups $H_1$ and $H_2$ in these conditions exist but are hard to compute, and so that the security of the subset version is greater than the subgroup version one even if the MSP answers \texttt{YES}. In some sense, having a difficult MSSP means that the subset version is harder than the subgroup one.
 
 \begin{problem}
 Study the MSP and the complexity of the MSSP for different classes of groups.
 \end{problem}
 
 However, the standard method for finding commuting subgroups in braid groups seems to be taking generators apart from each other. This is not an ideal strategy for us, unless we can somehow disguise the two disjoint subsets of generators used to construct  the algebraic subsets $S$ and $T$. In the Thompson's group $F$ the strategy is similar.

\subsubsection{Second protocol}

 In this case, the same objection would apply if one could find subgroups $H_1,H_2\leq G$ such that $S\subseteq H_1 \subseteq C(a_1)$ and $T\subseteq H_2 \subseteq C(b_2)$. However, since $a_1$ and $b_2$ are not known by the adversary Eve, it is not clear what subgroups Eve should look for and for this reason, the strengthening of the protocol by taking subsets is clearer in this version. However, since $S\subseteq C(a_1)$, then $\langle S\rangle \leq C(a_1)$, so it is important to ensure that $\langle S\rangle$ is not finitely generated, or, in case it is, that it cannot be easily computed. Again, by what was seen above, $\langle S\rangle$ is not finitely generated in general and the problem of determining if it is finitely generated or not seems to be unexplored in group theory.
 \begin{problem}
 For a group $G$, can we describe which algebraic subsets $S$ are such that $\langle S \rangle$ is finitely generated? Is it decidable, on input $S\in \Alg(G)$, whether $\langle S\rangle$ is finitely generated or not and, in case it is, can we compute it?
 \end{problem}
 
 \subsection{Some examples}
 The goal of this subsection is to highlight the potential of this approach in comparison with the traditional one. We do not make claims about the security of the protocol in these cases nor propose these groups as a suitable platform, but highlight the range of possibilities provided by our approach. 
  \subsubsection{Ascending HNN-extensions} \label{sec: ahnn}
  Consider the Baumslag--Solitar group $BS(1,m)$, $m\geq 1$, defined by the following presentation:
  \begin{align*}
BS(1,m)=\langle\, a,t \mid t^{-1}at= a^m\rangle.
\end{align*}

These are precisely the metabelian (and the solvable) Baumslag--Solitar groups. Let $N$ be the commutator subgroup of $BS(1,m)$. Then, for any element $x$ in $BS(1,m)$, we have that the centralizer
\begin{align*}
C(x)=\begin{cases}
N \quad &\text{if $x\in N$}\\
\text{cyclic} & \text{if $x\not\in N$},
\end{cases}
\end{align*}
and $N$ is isomorphic to the additive group of $m$-adic rational numbers $\Z[\frac 1 m]$, which is not finitely generated (see, for example, \cite{[CK12]} for a reference).
However, they are always algebraic.

\begin{proposition}
Let $m\geq 1$ and $G=BS(1,m)$. Then, the centralizer of any element in $G$ is algebraic.
\end{proposition}
\noindent\textit{Proof.}
Let $x\in G$ and $N$ denote the commutator subgroup of $G$. If $x\not\in N$, then $C(x)$ is cyclic, and so it is rational (thus algebraic). Let $S= \{t^{-k}at^k \mid k\in \Z\}$. It is well known that $S\in \Alg(G)$. From Lemma \ref{lem: generated algebraic}, it suffices to see  that  $N=\langle S\rangle$. Standard techniques show that  $N$ is the kernel of the natural homomorphism onto the infinite cyclic group obtained by setting $a=1$ in the presentation, and so it is the normal closure of $a$ (see \cite{[Co01]}), which is precisely $\langle  S\rangle$. 
\qed\\

However, despite illustrating the greater generality of the class of algebraic subsets when compared to the class of finitely generated subgroups, this does not serve our purpose, since $N$ is normal and $N$ is the centralizer of any element in $N$. So, the decomposition search problem would reduce to, given $u,w\in G$, find $a,b\in N$ such that $awb=u$, provided that such a pair of elements exists. But this equality can be rewritten as $w^{-1}awb=w^{-1}u$ and, since $N$ is normal, $w^{-1}aw$ and $b$ both belong to $N$, so also does $w^{-1}u$. Hence, we could simply take $a=1$ and $b=w^{-1}u$.

This is part of a more general setting, as the metabelian Baumslag--Solitar groups are ascending HNN-extensions of $\Z$.
  Let $G$ be a group and $\phi$ be an injective endomorphism of $G$. 
  The \emph{ascending HNN-extension of $G=\langle X\mid R\rangle$ induced by $\phi$} is the group with the following presentation:
  \begin{align*}
G\ast_\phi =\left\langle X,t \left\vert \begin{array}{ll} R&
 \\ t^{-1}xt= x\phi, & x\in G
\end{array} \right. \right\rangle
\end{align*}
   
In \cite{[Car23c]}, the first author remarks that all orbits of elements by automorphisms of a group $G$ are algebraic subsets of a cyclic extension of $G$, using this to prove that a kind of Fatou property does not hold for algebraic subsets, contrary to what happens for rational subsets. The same idea can be applied for injective endomorphisms and ascending HNN-extensions.
\begin{proposition}
Let $G$ be a group, $x\in G$ and $\phi$ be an injective endomorphism of $G$. Then the orbit of $x$ through $\phi$, $\Orb_\phi(x)=\{x\phi^k\mid k\in \N\}$ is an algebraic subset of $G\ast_\phi$.
\end{proposition}
\noindent\textit{Proof.}  Let $A$ be a generating set of $G$. By the fact that $t^{-1}xt= x\phi$, we deduce that $t^{-k}xt^k= x\phi^k$, for all $k\in \N$, and so $\{t^{-k}xt^k\mid k\in \N\}\pi=\Orb_\phi(x)$. Since $\{t^{-k}xt^k\mid k\in \N\}\subseteq \widetilde{(A\cup t)}^*$ is a context-free language, then $\Orb_\phi(x)\in\Alg(G\ast_\phi).$
\qed\\

 \subsubsection{Hiding the generators}
 The braid group on $n$ strands is given by the following presentation
\begin{align*}
B_n =\left\langle \sigma_1,\sigma_2,\ldots ,\sigma_{n-1} \left\vert \begin{array}{ll} \sigma_i \sigma_j =\sigma_j
\sigma_i & |i-j|> 1 \\ \sigma_i \sigma_{j} \sigma_i =\sigma_{j} \sigma_i \sigma_{j} & |i-j|=1
\end{array} \right. \right\rangle
\end{align*}
 As mentioned above, typical constructions of commuting subgroups in braid groups involve choosing generators apart from each other. If the same strategy is taken when considering subsets, we will end with two subsets $S\subseteq \langle\sigma_1,\ldots, \sigma_{i-1} \rangle$ and $T \subseteq \langle\sigma_{i+1},\ldots, \sigma_{n-1} \rangle$ contained in two pairwise commuting subgroups and the problem is reduced to the subgroup case. A possible strategy is to mask the generators while designing the context-free grammar. We remark that this is a particular case of the MSSP, which might be worth considering given the importance that braid groups have in protocols of this nature.
 
 The idea is that, to describe $S=L\pi$, rules can be added involving generators from $\{\sigma_{i+1},\ldots, \sigma_{n-1} \}$, making sure that they will end up cancelling (not as words, but as group elements), meaning that any derived word $w\in L$ is such that $w\pi \in  \langle\sigma_1,\ldots, \sigma_{i-1} \rangle$, but words in $\widetilde{\{\sigma_{i+1},\ldots, \sigma_{n-1} \}}^*$ appear in the grammar.
 The following problem is thus natural: 
 \begin{problem}
 Given two algebraic subsets $S,T\in \Alg(B_n)$ such that  $S\subseteq \langle\sigma_1,\ldots, \sigma_{i-1} \rangle$ and $T \subseteq \langle\sigma_{i+1},\ldots, \sigma_{n-1} \rangle$ for some $i\in\{2,\ldots,n-2\}$, what is the complexity of finding such an $i$?
\end{problem}

Again, we remark that this is not the problem which the protocol is based on, it just encodes the difficulty of reducing the problem to the subgroup case.

 We now present a possible implementation of this method by showing an example of a group with  pairwise commuting infinitely generated algebraic subgroups.

 \section{A protocol in ascending HNN-extensions of free-abelian groups}\label{proposal}
 
Let $H=\Z^m$, $u,v\in \Z^m$ be any vectors and $M$ be an integral $m\times m$ matrix with nonzero determinant and write $G=H\ast_M$. Put $S=\langle \{t^{-k}u t^k\mid k\in \Z\} \rangle$ and $T=\langle \{t^{-k}v t^k\mid k\in \Z\} \rangle$. By Lemma  \ref{lem: generated algebraic}, both these subsets are algebraic. We remark that groups (and subgroups) of this kind were studied in \cite{[Val19]}.
We will now see that they are not necessarily finitely generated and that they commute pairwise.
It is enough to check commutativity on the generators.  
We have that:
\begin{align}\label{c1}
(t^{-p}u t^p)(t^{-q}vt^q)=t^{-p}ut^{p-q}vt^{q} 
\end{align}
and 
\begin{align}\label{c2} 
(t^{-q}vt^q)(t^{-p}ut^p)=t^{-q}vt^{q-p}ut^{p}.
\end{align}

If $p=q=0$, then clearly (\ref{c1})=(\ref{c2}). If $p,q\geq 0$, then $t^{-p}ut^p=uM^p\in \Z^m$ and $t^{-q}vt^q=vM^q\in \Z^m$, and so they commute.

 If  $p,q\leq 0$ then, assuming w.l.o.g. that $p\geq q$,
 $$(\ref{c1})=t^{-p}t^{p-q}(uM^{p-q}) vt^{q} =t^{-q}(uM^{p-q})vt^q$$
and 
 $$(\ref{c2})=t^{-q}v(uM^{p-q})t^{q-p}t^{p}=t^{-q}(uM^{p-q})vt^q,$$
 so (\ref{c1})=(\ref{c2}). 
If $pq<0$, then, assume w.l.o.g. that $p>0$ and $q<0$. Then
$$(\ref{c1})=t^{-p}t^{p-q}(uM^{p-q})vt^{q}=t^{-q}(uM^{p-q})vt^q$$
and 
 $$(\ref{c2})=t^{-q}v(uM^{p-q})t^{q-p}t^{p}=t^{-q}(uM^{p-q})vt^q$$
 so (\ref{c1})=(\ref{c2}) and $S$ and $T$ are pairwise commutative.

Notice that, in the particular case where $M$ is a diagonal matrix, $u=(1,0,\ldots, 0)$ and $v=(0,1,0,\ldots, 0)$, then  $\langle u,t \rangle\simeq BS(1,M_{11})$ and $\langle v,t \rangle\simeq BS(1,M_{22})$ via the obvious isomorphisms and that $S$ and $T$ are the corresponding commutator subgroups, hence not finitely generated.

So, for an implementation, an injective matrix $M$ and elements $u,v\in \Z^m$ would be chosen and made public. Then, we construct context-free grammars representing the subsets $S$ and $T$, which are also made public together with an element $w\in G$. Alice would then pick $a_1\in S$ and $b_2\in T$, by choosing words derived using the respective context-free grammars. The product  $a_1wb_1$ would then be sent to Bob; proceeding similarly, Bob would pick $a_2\in S$ and $b_2\in T$ and send $b_2wa_2$ to Alice; both of them could now compute the shared secret key $a_1b_2wa_2b_1$.

 In the above examples, the subsets are also subgroups. This needs not to be the case. Obviously, since we require pairwise commutativity between the subsets, there must be some structure involved, namely concerning centralizers. We considered subgroups because they turned out to be algebraic, but in general, any algebraic subset of an adequate centralizer would suffice and this can be chosen in a language-theoretic way in order to decrease the level of structure in our subset.

\section{Length and  distance based attacks}\label{sec: attacks}

Most length and distance based attacks (see \cite[Section 11.2]{[MSU11]}) have the following general setting: there is a length (or distance) function $\ell$ and we start multiplying generators that minimize $\ell$ in order to obtain the shared secret key. Concretely, in \cite{[RST07]}, Ruinskiy, Shamir and Tsaban describe an attack to the  protocol proposed in \cite{[SU05]} (recall the first protocol described in this paper) by defining a function that, \emph{in some sense}, measures how close an element is to $B$ and then proceeding as follows: let $S_A=\{\alpha_1,\ldots, \alpha_k\}$ denote the generating set of $A$
\begin{enumerate}
\item Let $\tilde a=1$.
\item For every generator $\alpha_i\in S_A^{\pm 1}$, compute $a^{(i)}=\tilde a\alpha_i$ and $b^{(i)}=w^{-1}(a^{(i)})^{-1}a_1wb_1$, and measure how far $b^{(i)}$ is from $B$. If $b^{(i)}\in B$, then we can stop because we have found $a^{(i)}\in A$ and $b^{(i)}\in B$ such that $a^{(i)}wb^{(i)}=a_1wb_1$.
\item We choose $j$ to minimize the distance between $b^{(j)}$ and $B$, put $\tilde a=a^{(j)}$ and return to step 2.
\end{enumerate}
Similar approaches are taken to tackle different problems, like the (simultaneous) conjugacy search problem. 

We believe that our approach is more resistant to this kind of attacks than the subgroup one.  The main reason for this is that there is no clear way of enumerating the elements of the subset, as they might represent a non-finitely generated subgroup, or not even a subgroup, and multiplying by generators of the group might not work because we might end with a word not belonging to the subset $A$, and so, not commuting with $B$.

We are not aware of a similar strategy that can be taken in the context of algebraic subsets.
\begin{problem}
Can we develop length or distance based attacks for protocols relying on pairwise commuting algebraic subsets?
\end{problem}
The natural approach seems to be going down the derivation tree making the choice of path minimizing $\ell$ at each step, but it is not clear that it should work in the same way due to depth differences.
\section{More group-theoretic problems arising from this work}\label{sec: other problems}
In \cite{[Shp18]}, Shpilrain presents some problems in group theory motivated by cryptography. Interestingly, some nontrivial (and natural) problems like the factorization decision problem only arises in group theory with cryptography. Apart from all the problems discussed above, we present some new problems in group theory that are \emph{in some way} related to the discussion above.

In \cite{[Car23c]}, the first author remarks that all orbits of elements by automorphisms of a group $G$ are algebraic subsets of a cyclic extension of $G$. The idea is the same as the one in Section \ref{sec: ahnn}. Orbits by automorphisms and endomorphisms are of particular interest in group theory and orbit problems, inspired by the work of Brinkmann \cite{[Bri10]}, have deserved a lot of attention by researchers in group theory in the past years \cite{[KL86],[Bri10], [BMMV06],[BMV10], [SW10],[Ven14],[Car23], [Car23b],[Car23c],[Car23z],[Log22],[CD23a]}. 

So, we propose the following problems that arise from the conjugacy, decomposition and factorization problems restricted to orbits of elements for groups in a certain class (e.g. let $G$ be a virtually free group):
\begin{problem}[Conjugacy problem restricted to orbits]
Is there an algorithm taking as input an endomorphism $\phi\in \End(G)$ and elements $x,y,z\in G$ that decides whether $x$ and $y$ are conjugate by an element in the orbit of $z$, i.e., whether there is some $k\in\N$ such that $(z^{-1}\phi^k)x (z\phi^k)=y$?
\end{problem}
\begin{problem}[Decomposition problem restricted to orbits]
Is there an algorithm taking as input an endomorphism $\phi\in \End(G)$ and elements $x,y,z,w\in G$ that decides whether there are  $k_1,k_2\in\N$ such that $(z\phi^{k_1})x (w\phi^{k_2})=y$?
\end{problem}
\begin{problem}[Factorization problem restricted to orbits]
Is there an algorithm taking as input an endomorphism $\phi\in \End(G)$ and elements $x,y,z\in G$ that decides whether there are  $k_1,k_2\in\N$ such that $(y\phi^{k_1})(z\phi^{k_2})=x$?
\end{problem}

\begin{problem}
How difficult are the search versions of these problems?
\end{problem}

Computing high powers of endomorphisms is usually computationally demanding. However, in \cite{[HKKS13]}, a cryptosystem based on semidirect products and involving the computation of powers of automorphisms is proposed. If orbit search problems are hard, then  one could consider a Diffie--Hellman-like key exchange protocol based on the search version of Brinkmman's equality problem, which consists on deciding, given $\phi\in \End(G)$ and $x,y\in G$, whether there is some $k\in \N$ such that $x\phi^k=n$:

\begin{center}
\fbox{\begin{minipage}{\textwidth}
\noindent\textbf{Public information:} group $G$, $x\in G$, $\varphi\in \End(G)$. \\

\noindent\textbf{Key exchange:}
\begin{enumerate}
\item Alice picks $m\in \N$, computes $x\phi^m$ and sends it to Bob.
\item Bob picks $n\in \N$, computes $x\phi^n$ and sends it to Alice.
\item Alice computes  $x\phi^n\phi^m=x\phi^{n+m}$ and Bob computes  $x\phi^m\phi^n=x\phi^{n+m}$. 
\end{enumerate}

\noindent\textbf{Shared secret key:}  $K=x\phi^{n+m}$.\\
\end{minipage}}
\end{center}

\begin{problem}
What properties must the group $G$ and the endomorphism $\phi$ have to make the above protocol secure?
\end{problem}
\section{Conclusion}\label{sec: conclusion}
We propose the introduction of subsets in some public-key cryptosystems. This seems natural, as decision problems are in general harder for the subset variants of the problems, so it makes sense to consider the search versions as well. We consider algebraic subsets because they seem to be the most balanced in terms of simplicity of description vs. difficulty, but, in theory,  one could go beyond that in Chomsky's hierarchy and consider, for instance, recursively enumerable subsets. We highlight in Section \ref{sec: attacks} that these subsets seem to be more resistant to traditional length-based attacks. This paper contains more questions than answers and has the goal of promoting the discussion around the security of these variations of known protocols. Additionally, the orbit problems proposed in Section \ref{sec: other problems} should be interesting to the group theoretical community.

\bibliographystyle{plain}
\bibliography{Bibliografia}

 \end{document}